# A possible use of the cubit rod in Ancient Egypt to measure and draw lengths based on fractions.

*Massimiliano Benes*

*Università degli Studi di Udine*    massimiliano.benes@uniud.it

*Abstract*

*We will discuss about a possible method of using the cubit rod by the architects and the surveyors of Ancient Egypt to measure and draw lengths, comparing it with the other interpretations present in Literature. Instead of the modern decimal notation, which sees the use of comma to represent a number or a measure, at that time there was a wide use of fractions in calculations. The current work proposes that, through the cubit rod and its partitions of the finger into fractions, it could be possible to obtain very accurate measurements.*

*Introduction*

In the Ancient Egypt the main unit of measurement was the cubit, based on the forearm length from the tip of the medius finger to the bottom of the elbow, equivalent to about 52,5 cm.
The royal cubit, also known as pharaonic, was divided into 28 fingers (the finger was 18,75 mm wide) but has other partitions, for example the palm corresponded to 4 fingers and the fist to 6 fingers. The short cubit was instead divided into 24 fingers [1, 2]. This unit of measurement was made materially by creating real rulers, the so-called cubit rods, which could also be foldable.
At the archaeological excavations several cubit rods in stone, metal or, more rarely in wood covered with gold, were rediscovered.

One of the rediscovered wood cubit rod is the one in the grave of the architect Kha, who lived between 1400 and 1350 BC. which was found by Ernesto Schiapparelli at the archaeological site located north of Deir el-Medina in 1906 and exhibited at the Egyptian museum in Turin, it probably had a more ceremonial than practical significance.

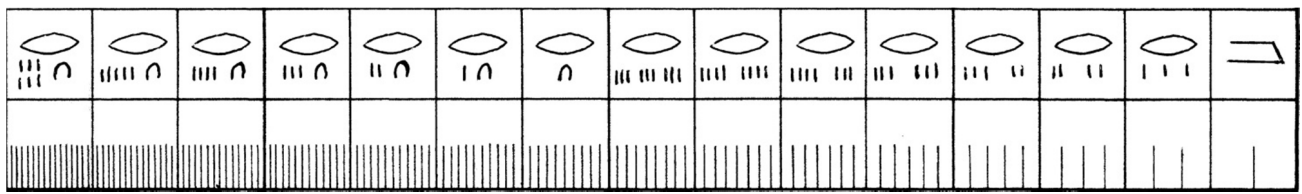

Figure 1: reproduction of a cubit rod vertical face detail.

The figure 1 represents a reproduction of a cubit rod vertical face detail. Looking at the figure (first line at the bottom) it is possible to see the partition of the scale in fingers, but it is also possible to observe a finer partition: in fractions of fingers for the first 15. The first finger, the one on the right is divided into 2 equal parts, the second in 3, the third in four and so on up to the fifteenth, which is divided into 16 equal parts, evidently the maximum resolution allowed or required at that time.

In the second line starting from the bottom the corresponding numeration is shown, according to the symbols in use.

Based on the illustration, one wonders which could have been the use of a ruler like that and why the partition of the fingers started from 2 without a first whole finger appearing in the cubit rod.

*Description of the measurement method using fractions*

Probably the ancient Egyptians measured lengths based on the fractions by an appropriate reading of the ruler scale, as presented below.
We must think that ancient Egyptian, even though they knew the base ten and therefore the exponents of ten to count quantities higher than 1, they used to resort to fractions to deal with quantities lower than 1.
Fractions were used in place of the current decimal numeration consisting in using the comma to represent a number or a measure.
For example, thanks to the papyrus of Rhind (one of the rarest rediscovered papyrus containing mathematic topics) [3], the sixth problem of the partition of 9 loaves among 10 workers is known.
This problem was solved by ancient Egyptian adding up fractions: the ration of bread getting to each worker was given by: 2/3+ 1/5+ 1/30 of loaf.
In this papyrus only algebra is used to solve the problem employing, in contrast to the major possibilities of cubit rod, only unitary fractions, save 2/3 [4].

Before facing up with the measurement problem, should deal with a preliminary aspect: it is possible that the lack of the first whole finger on the cubit rod scale could be linked to the absence of the concept of zero in Egyptian culture.
Nowadays the measurement of a segment using a ruler involves the employment of the concept of zero when overlapping the beginning of the scale, the zero indeed, with the beginning of the segment, starting from that reference point to count notches.
In absence of the concept of zero, there cannot be used a procedure like that in measuring the length of a segment or an object through the ancient ruler. Cubit rods having a specular partition of the scale compared to figure 1, also were found, as for example the royal cubit rod of Kha, but as we will see later this is not going to change our discussion.

It comes to think that ancient Egyptian understood the concept of zero but they did not push in dividing the cubit beyond 1/16 of finger (about 1,2 mm), even though the presence of studies, based on historical sources, which indicate the knowledge of inferior fractions such as 1/32 and 1/64 or even lower [3]. Ancient Egyptian could have approached the concept of zero through a progressive reasoning typical of the limit, by dividing a range in even smaller parts.
Ancient Egyptian had already faced up with this topic in adding the fractions like $1/2^n$ (with n greater than 0), which has the unit as limit and the reasoning is represented with the graphic construction of the "Eye of Horus" [5].

Let's resume the initial problem about how to measure a length with the cubit rod by using fractions and without the concept of zero.
For convenience's sake we will call "notch" the line engraved on the ruler which separates a finger from another and "incision" each line dividing the finger into its fractions.

Therefore, supposing to have to measure a segment of a certain unknown length, it is always possible to proceed manually with the ruler until we find the best correspondence of a notch with an extreme of the segment and of an incision with the other extreme of the segment.

This can be done proceeding with both the alignments, starting from the left (dotted line) or from the right (solid line) and starting from whichever notch on the ruler scale.
In order to better illustrate the proposed procedure, an example is shown in figure 2.

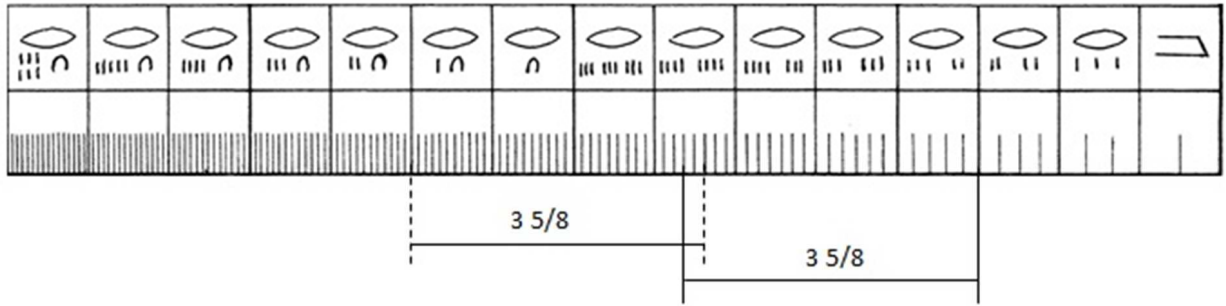

Figure 2: Example of measurement of the length of a segment with the cubit rod.

From the observation of the previous figure, a measurement of the segment corresponding to 3 whole fingers and 5/8 of finger is obtained as there is a correspondence with the fifth incision of the eighth finger: with the ruler it was possible to read the denominator of the fraction directly referring to the symbols reported above the scale. In this case, since there is no prime number in the denominator, it is possible to write the result of the measurement in unitary fractions as 3 1/2 1/8 (according to the notation of that time).

In the same way the ruler could be used to draw a segment of any length that was known as the sum of a whole number and a fraction: an example is given in figure 3, where it is attempted to trace a segment 2 fingers and 4/5 of finger long. Therefore, had been considered 2 notches and 4 incisions.

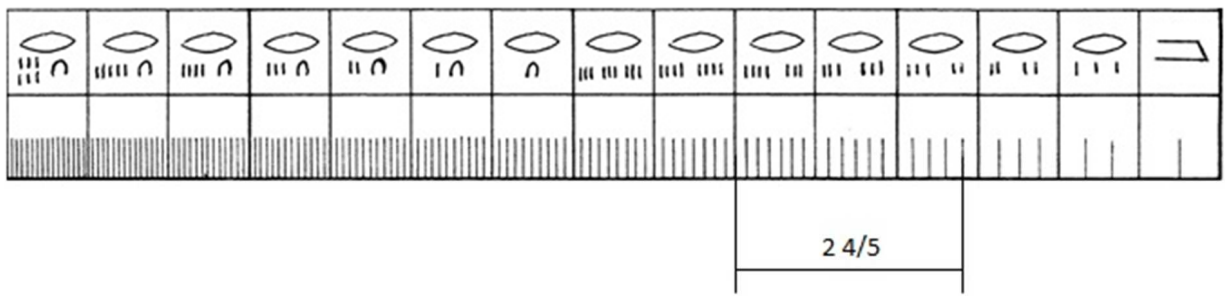

Figure 3: Example of tracing of the length of a segment with cubit rod.

Referring to figure 3, it is to be noted that in the second case we carry on from left to right as opposed to the first case, to demonstrate the equivalence of the way to proceed.
While, in the case of very long segments, it is possible to let the segment, to measure or trace, start from the not graduated scale, available on the left of the graduated one, unless resorting to a double measurement.

*Discussion*

The presented method is a proposal trying to interpret how measurements or drawing of segments could be conducted in ancient times, through the use of cubit rod so by using a scale divided into fractions instead of decimal parts.
As illustrated, the measurement takes place with the only ruler without the aid of further instruments or graduated bars. Even though it is possible that they were used, as supposed by [6], who illustrates the use of the ruler according to a method similar to the one of the caliper of Vernier for the measurement of the length of objects.

Of course, the more attempts were made and the greater was the experience of the architect or surveyor and the greater and greater precision could be obtained in measuring.
It is possible that these measurements were carried out on drawings and then communicated the measures to the operational teams, perhaps by employing an appropriate scaling coefficient: for instance, 3 fingers and 5/8 of finger could be the equivalent of 3 cubits and 5/8 of cubit in reality of construction.
It is to be noted that a message like that, in the form of fraction, could present fewer communication errors and could be more direct compared to the evaluation of the correct position of the comma, characteristic of the use of decimal digits, in case of a change in the unit of measurement.

In the previous example, since 3+ 5/8 corresponds to 3,625, the way of representation by fractions is much more immediate and precise than the second one, on which burden the problem of establishing a shared method about the significant digits to take into account, with consequences in the propagation of approximations for the realization of the works.
On the other hand, some limitations in the accuracy of measurement could derive from the level of precision with which the incisions on the rulers were made by the artisan of that time.
It is likely anyway, considering the remarkable architectural results obtained by ancient Egyptian, that the employed manufacturing precision for the production of cubit rods and their samples was more than sufficient.

To support our description of the method of measurement, the statues of the Sumerian architect king Gudea (2175 BC), brought back to lite by Ernest de Sarzec at the archeological site of Lagash in 1881, hold on their knees some small tablets with a bas relief, illustrating a graduated ruler and a drawing stylus [7, 8], demonstrating that the practise of tracing segments by using graduated rulers, dates back to previous times, prior to the Egyptian cubit rods.
This ruler is divided into 16 fingers with an average length of 16.8 mm [5] and only 5 fingers are further divided into fractions, up to 1/6, interchanging a divided finger with a not divided one.
The last finger also shows, on the opposite side of the scale, a further division of 1/6 into two and three parts, always inserting an empty space. The method of measurement which has been discussed is also applicable to the ruler of Gudea. In any event, ancient Egyptian have surpassed in precision and completeness the ruler of Gudea, which lacks the subsequent fractions to 1/6.

*Conclusions*

In the current work, a possible use of the cubit rod to measure and draw the length based on fractions, a mathematical concept well known by ancient Egyptians, has been explained.
This method is based on the use of the only ruler and allows to obtain very precise measure of length, without resorting to the concept of zero. The communication of a measure, read on a project drawing on the construction site, could be considerably facilitated by the fractions that automatically allowed the application of an appropriate scaling factor to the anthropometric unit of measurements of that time.

The ruler of Gudea seems to be based on the same logic, but the ancient Egyptian surpassed it in precision and completeness. Finally, although this does not emerge from the papyrus of Rhind, it is possible that the measurements of precision with the ruler could be useful to carry out geometry studies and to trace segments of a known length in the form of fraction. Conceptual instruments which were within the reach of the mathematician of that time and handled with great skill.